\documentclass[11pt,a4paper]{article}
\usepackage[T1]{fontenc}
\usepackage{lmodern}
\usepackage[margin=28mm]{geometry}
\usepackage{amsmath,amssymb,amsthm}
\usepackage{graphicx}
\usepackage{booktabs}
\usepackage{authblk}
\usepackage[hidelinks]{hyperref}

\newtheorem{theorem}{Theorem}[section]
\newtheorem{lemma}[theorem]{Lemma}
\newtheorem{proposition}[theorem]{Proposition}
\newtheorem{corollary}[theorem]{Corollary}
\theoremstyle{definition}
\newtheorem{construction}[theorem]{Construction}
\theoremstyle{remark}
\newtheorem{remark}[theorem]{Remark}

\newcommand{\Tr}{\mathcal{T}_r}
\newcommand{\Ep}{\mathcal{E}^{(r)}_{p,s}}
\newcommand{\M}{\mathcal{M}}
\hypersetup{pdftitle={The isolation number of uniform hypertrees}}

\begin{document}
\title{The isolation number of uniform hypertrees\thanks{This work is supported by the National Natural Science Foundation of China (No.~12461065).}
}

\author[1\thanks{Email: liminhuimath@163.com.}]
{Minhui Li}

\author[2,3\thanks{Corresponding author. Email: zhsm1@qhnu.edu.cn.}]
{Shumin Zhang}

\author[2\thanks{Email: qujunhaosx@163.com.}]
{Junhao Qu}

\affil[1]{\footnotesize Center for General Education, Qinghai Institute of Technology, Xining, Qinghai 810006, China}

\affil[2]{\footnotesize
School of Mathematics and Statistics, Qinghai Normal University, Xining, Qinghai 810001, China}

\affil[3]{\footnotesize
The State Key Laboratory of Tibetan Intelligent Information Processing and Application, Xining, Qinghai 810008, China}

\maketitle

\begin{abstract}
An isolating set of a hypergraph $H$ is a vertex set whose closed neighbourhood intersects every edge of $H$.
The isolation number $\iota(H)$ is the minimum cardinality of such a set.
Let $T$ be an $r$-uniform hypergraph whose incidence graph is a tree.
We prove that if $T$ has order $n$ and at least two edges, then
$\iota(T)\le\lfloor n/(2r-1)\rfloor$, where $r\ge3$.
Moreover, the bound is attained for every admissible order.
Furthermore, we characterize all hypertrees attaining it by a construction using six types of marked hypertrees of orders $2r-1$, $3r-2$ and $4r-3$.
For connected graphs of order $n$ with at least two blocks, all of which are triangles, we obtain the sharp bound $\lfloor n/5\rfloor$ on the cycle isolation number and characterize the graphs attaining it.
\end{abstract}


\section{Introduction}

Let $H=(V,E)$ be a finite hypergraph, where $E$ is a family of nonempty subsets of $V$.
Its \emph{order} and \emph{size} are $n(H)=|V|$ and $m(H)=|E|$, respectively.
An edge of size $r$ is an \emph{$r$-edge}, and $H$ is \emph{$r$-uniform} if every edge is an $r$-edge. The hypergraph is \emph{linear} if any two distinct edges have at most one common vertex.
Two distinct vertices are \emph{adjacent} if they lie in a common edge.
For $v\in V$, let $N_H(v)$ be its neighbourhood and $N_H[v]=N_H(v)\cup\{v\}$.
For $S\subseteq V$, define $N_H[S]=\bigcup_{v\in S}N_H[v]$.
We omit the subscript when the hypergraph is clear from the context.

For $X\subseteq V$, the \emph{induced subhypergraph} $H[X]$ has vertex set $X$ and edge set $\{e\in E:e\subseteq X\}$.
Thus deleting vertices does not replace an edge by its intersection with the remaining set.
A set of vertices is \emph{independent} if it contains no edge.
Following the isolation terminology of Li, Zhang and
Ye~\cite{LiZhangYe}, we call $S\subseteq V$ \emph{isolating} if
$V\setminus N_H[S]$ is independent.
Equivalently, $N_H[S]$ intersects every edge of $H$.
The minimum cardinality of an isolating set is the \emph{isolation number} of $H$, denoted by $\iota(H)$.
This condition is weaker than domination, which requires $N_H[S]=V$.

Caro and Hansberg~\cite{CaroHansberg} studied isolation in graphs
as a form of partial domination. Lema\'nska et al.~\cite{LemanskaTrees} characterized the trees $T$ satisfying $\iota(T)=n(T)/3$. More generally, a set $S$ is \emph{$\mathcal{F}$-isolating} if $H[V\setminus N_H[S]]$ contains no copy of any hypergraph in $\mathcal{F}$. Li, Zhang and Ye~\cite{LiZhangYe} developed this notion for hypergraphs.
Write $K_k^r$ for the complete $r$-uniform hypergraph on $k$ vertices,
and $\iota(H,K_k^r)$ for the minimum size of a $\{K_k^r\}$-isolating set.
In a recent preprint, Borg~\cite{BorgCliques} proved that for any connected $n$-vertex $r$-uniform hypergraph,
\[
\iota(H,K_k^r)\le\frac{n}{k+1}\qquad(2\le r\le k),
\]
with two exceptions: $H\cong K_k^r$, and, in the graph case $k=r=2$, $H\cong C_5$.
Taking $k=r\ge3$ gives $\iota(H)\le n/(r+1)$ unless $H$ consists of a single edge. Moreover, building on known results for the graph case $r=2$,
the paper determines the extremal structures for general $r$~\cite{BorgCliques}.
More recently, Baluci and Borg~\cite{BaluciBorg}
proved a sharp upper bound for the isolation number
of connected uniform hypergraphs in terms of their size.
In particular, taking $k=r$ in their result gives
\[
\iota(H)\le
\left\lfloor\frac{(r-1)m(H)+1}{2r-1}\right\rfloor
\]
for every connected $r$-uniform hypergraph $H$
with $H\not\cong K_r^r$.
By~\eqref{eq:order}, this bound yields~\eqref{eq:bound}
for every $T\in\Tr$. Our main objective is to characterize all hypertrees
attaining equality in~\eqref{eq:bound}.

A graph is a \emph{forest} if it contains no cycle, and a connected forest is a \emph{tree}. We allow a tree to consist of a single vertex.
The \emph{incidence graph} $B(H)$ is the bipartite graph with vertex classes $V(H)$ and $E(H)$, where $v$ and $e$ are adjacent if and only if $v\in e$.
We call $H$ a \emph{hypertree} if $B(H)$ is a tree.
For the uniform hypergraphs considered here, this definition
is equivalent to Bacher's notion of an acyclic
hypertree~\cite[Section~2]{Bacher}; see also Masbaum and
Vaintrob~\cite[Definition~3.2]{MasbaumVaintrob}
for the $3$-uniform case.
We specify this convention because the term ``hypertree''
has other meanings in the literature.

Equivalently, a hypertree is connected and has no Berge cycle, where a \emph{Berge cycle} of length $k\ge2$ consists of distinct vertices $v_1,\ldots,v_k$ and distinct edges $e_1,\ldots,e_k$ with $\{v_i,v_{i+1}\}\subseteq e_i$ for every $i$, taking $v_{k+1}=v_1$.
In particular, every hypertree is linear: two distinct edges sharing
two vertices would give a $4$-cycle in its incidence graph.

Throughout the paper, $\Tr$ denotes the family of $r$-uniform hypertrees with at least two edges, where $r\ge3$.
For $T\in\Tr$, $B(T)$ is a tree with
$n(T)+m(T)$ vertices and $rm(T)$ edges.
Thus $rm(T)=n(T)+m(T)-1$, which gives
\begin{equation}\label{eq:order}
n(T)=(r-1)m(T)+1\qquad(T\in\Tr).
\end{equation}

The result of Baluci and Borg~\cite{BaluciBorg}
gives the following bound for $r$-uniform hypertrees.
We give a direct proof based on edge-neighbourhood
packings, which will also be used in the extremal
characterization.
\begin{theorem}\label{thm:bound}
Let $T\in\Tr$ be a hypertree of order $n$. Then
\begin{equation}\label{eq:bound}
\iota(T)\le\left\lfloor\frac{n}{2r-1}\right\rfloor.
\end{equation}
For any $n\ge2r-1$ with $n\equiv1\pmod{r-1}$, there is a member of $\Tr$ of order $n$ attaining the bound.
\end{theorem}

Theorem~\ref{thm:extremal} characterizes all hypertrees attaining
equality in~\eqref{eq:bound} by a construction using marked
hypertrees and an auxiliary hypertree. The construction uses
six types of marked hypertrees for $r\ge3$.
Corollary~\ref{cor:unrounded} characterizes the
hypertrees satisfying $\iota(T)=n(T)/(2r-1)$.
We also apply these results to cycle isolation in
connected graphs whose blocks are all triangles.

The proof uses the classical transversal--packing property of subtrees of a tree. We give the required argument for completeness; more general subtree packing and covering problems were studied by B\'ar\'any, Edmonds and Wolsey~\cite{BaranyEdmondsWolsey}.
The additional structure needed for the extremal characterization comes from the edge neighbourhoods of a uniform hypertree.

\section{Edge neighbourhoods and the isolation number}
We first recall a property of subtrees of a tree.
Let $\mathcal{F}$ be a finite family of nonempty subtrees of a
tree $G$. The \emph{transversal number} of $\mathcal F$ is defined by
\[
\tau(\mathcal F)
=
\min\bigl\{
|S|:S\subseteq V(G),\
S\cap V(Q)\ne\emptyset
\text{ for every }Q\in\mathcal F
\bigr\}.
\]
The \emph{packing number} of $\mathcal F$ is defined by
\[
\nu(\mathcal F)
=
\max\bigl\{
|\mathcal P|:\mathcal P\subseteq\mathcal F,\
V(Q)\cap V(R)=\emptyset
\text{ for all distinct }Q,R\in\mathcal P
\bigr\}.
\]
For vertices $u,v\in V(G)$, let $d_G(u,v)$ denote
the distance between $u$ and $v$ in $G$, that is,
the length of a shortest $u$--$v$ path.
\begin{lemma}\label{lem:subtrees}
Let $\mathcal F$ be a finite family of nonempty subtrees of
a tree $G$. Then $\tau(\mathcal F)=\nu(\mathcal F)$.
\end{lemma}
\begin{proof}
If $\mathcal F=\emptyset$, then
$\tau(\mathcal F)=\nu(\mathcal F)=0$.
Assume that $\mathcal F\ne\emptyset$.

Choose an arbitrary vertex $v_0$ of $G$ as the root.
For each $R\in\mathcal F$, let $v_R$ be the unique vertex of
$R$ at minimum distance from $v_0$.
Let $Q\in\mathcal F$ be such that $d_G(v_0,v_Q)$ is maximum.

Let $R\in\mathcal{F}$ satisfy $V(R)\cap V(Q)\ne\emptyset$, and $u\in V(R)\cap V(Q)$.
Both $v_R$ and $v_Q$ lie on the $v_0$--$u$ path.
Since $d_G(v_0,v_R)\le d_G(v_0,v_Q)$, the $v_R$--$u$ path
contains $v_Q$. This path lies in $R$, so $v_Q\in V(R)$.
Thus every member of $\mathcal F$ intersecting $Q$ contains $v_Q$.

Select $v_Q$ and $Q$, remove all subtrees containing $v_Q$,
and repeat the same procedure on the remaining family until
no subtrees remain.
Let $Q_1,\ldots,Q_k$ be the subtrees selected at successive
steps, and let $S=\{v_{Q_1},\ldots,v_{Q_k}\}$ and $\mathcal P=\{Q_1,\ldots,Q_k\}$.
Every member of $\mathcal F$ contains a selected vertex,
so $S$ is a transversal of $\mathcal F$.
At each step, every subtree in the current family intersecting
the chosen subtree is removed. Hence the members of
$\mathcal P$ are pairwise vertex-disjoint.
Since $v_{Q_i}\in V(Q_i)$ for any $i$, the vertices
$v_{Q_1},\ldots,v_{Q_k}$ are distinct, and thus
$|S|=|\mathcal P|=k$.

Every transversal contains at least one vertex from each member
of any pairwise vertex-disjoint subfamily. Therefore,
\[
k=|\mathcal P|
\le \nu(\mathcal F)
\le \tau(\mathcal F)
\le |S|=k.
\]
Thus $\tau(\mathcal F)=\nu(\mathcal F)$.

\end{proof}

For an edge $e$ of a hypertree $T$, let $Q_e$ be the subtree of
$B(T)$ induced by the vertices at distance at most $3$ from
the vertex corresponding to $e$.
We call a set $P\subseteq E(T)$ an \emph{edge-neighbourhood packing}
if $N_T[e]\cap N_T[f]=\emptyset$ for all distinct $e,f\in P$.
\begin{theorem}\label{thm:packing}
For every $T\in\Tr$,
\[\iota(T)
=
\max\bigl\{|P|:
P\text{ is an edge-neighbourhood packing of }T\bigr\}.
\]
\end{theorem}

\begin{proof}
Let $\mathcal F=\{Q_e:e\in E(T)\}$.
We use the same notation for an edge of $T$ and its
corresponding vertex in $B(T)$.
For $v\in V(T)$ and $e\in E(T)$,
$N_T[v]\cap e\ne\emptyset$ if and only if $v\in N_T[e]$.
By the definition of $Q_e$, $N_T[e]=V(Q_e)\cap V(T)$.
Hence, a set $S\subseteq V(T)$ is an isolating set of
$T$ if and only if it is a transversal of $\mathcal F$.
Thus $\tau(\mathcal F)\le\iota(T)$.

Let $X$ be a minimum transversal of $\mathcal F$.
For each edge $f\in E(T)$ whose corresponding vertex
in $B(T)$ belongs to $X$, choose a vertex $u_f\in f$.
Let $S=(X\cap V(T))\cup\{u_f:f\in X\cap E(T)\}$.
If $f\in X\cap E(T)$ belongs to $V(Q_e)$, then
$d_{B(T)}(f,e)\le2$, since this distance is even and at most $3$.
Hence
\[
d_{B(T)}(u_f,e)
\le 1+d_{B(T)}(f,e)
\le3,
\]
so $u_f\in V(Q_e)$.
Thus $S$ is a transversal of $\mathcal F$ contained in $V(T)$.
Therefore, $S$ is an isolating set of $T$.
It follows that
\[
\iota(T)\le |S|\le |X|=\tau(\mathcal F).
\]
Hence $\iota(T)=\tau(\mathcal F)$.

We next show that for any $e,f\in E(T)$, $V(Q_e)\cap V(Q_f)\ne\emptyset$
if and only if $N_T[e]\cap N_T[f]\ne\emptyset$.
If $N_T[e]\cap N_T[f]\ne\emptyset$, then
$Q_e$ and $Q_f$ intersect, since
$N_T[e]=V(Q_e)\cap V(T)$ and
$N_T[f]=V(Q_f)\cap V(T)$.
For the converse, suppose that
$V(Q_e)\cap V(Q_f)\ne\emptyset$. Let $x\in V(Q_e)\cap V(Q_f)$.
If $x\in V(T)$, then $x\in N_T[e]\cap N_T[f]$.
If $x\in E(T)$, then both $d_{B(T)}(x,e)$ and $d_{B(T)}(x,f)$ are even.
Since both are at most $3$, they are at most $2$.
For any $u\in x$, we have $d_{B(T)}(u,e)\le3$ and $d_{B(T)}(u,f)\le3$.
Thus $u\in V(Q_e)\cap V(Q_f)$, so $u\in N_T[e]\cap N_T[f]$.

It follows that $\nu(\mathcal F)$ is the maximum cardinality
of an edge-neighbourhood packing of $T$.
By Lemma~\ref{lem:subtrees},
\[
\iota(T)=\tau(\mathcal F)=\nu(\mathcal F).
\]
\end{proof}

A set $P\subseteq E(T)$ is an edge-neighbourhood packing
if and only if $d_{B(T)}(e,f)\ge 8$
for any distinct $e,f\in P$.
The proofs of Lemma~\ref{lem:subtrees} and Theorem~\ref{thm:packing} yield a minimum isolating set and an edge-neighbourhood packing of the same cardinality.

For $e\in E(T)$, let $a_T(e)=|\{f\in E(T)\setminus\{e\}:f\cap e\ne\emptyset\}|$.
We abbreviate $a_T(e)$ as $a(e)$ when $T$ is clear from context.

\begin{lemma}\label{lem:neighbourhood}
Let $T\in\Tr$ and $e\in E(T)$. Then
$E(T[N_T[e]])=\{f\in E(T):f\cap e\ne\emptyset\}$
and
\begin{equation}\label{eq:neighbourhood}
|N_T[e]|=r+(r-1)a(e)\ge 2r-1.
\end{equation}
Moreover, $|f\cap N_T[e]|\le 1$ for any
$f\in E(T)\setminus E(T[N_T[e]])$.
\end{lemma}

\begin{proof}
Let $\mathcal E=\{f\in E(T):f\cap e\ne\emptyset\}$.
Since $B(T)$ is a tree, $|f\cap e|=1$ for any
$f\in\mathcal E\setminus\{e\}$.
Moreover, $(f\setminus e)\cap(g\setminus e)=\emptyset$
for distinct $f,g\in\mathcal E\setminus\{e\}$;
otherwise, $B(T)$ would contain a cycle.
By the definition of $N_T[e]$,
$N_T[e]=e\cup\bigcup_{f\in\mathcal E\setminus\{e\}}(f\setminus e)$.
Therefore,
\[
|N_T[e]|
=
r+\sum_{f\in\mathcal E\setminus\{e\}}|f\setminus e|
=
r+(r-1)a(e).
\]
Since $T$ is connected and has at least two edges,
$a(e)\ge 1$, which proves~\eqref{eq:neighbourhood}.

Let $B_e$ be the subgraph of $B(T)$ induced by the
vertices in $N_T[e]$ and the vertices corresponding
to the edges in $\mathcal E$.
The graph $B_e$ is connected, so $B_e$ is a tree.
Let $f\in E(T)\setminus\mathcal E$.
If $|f\cap N_T[e]|\ge 2$, then choose distinct
$x,y\in f\cap N_T[e]$.
The $x$--$y$ path in $B_e$, together with the path
$x,f,y$ in $B(T)$, forms a cycle in $B(T)$, a contradiction.
Thus $|f\cap N_T[e]|\le 1$ for any $f\in E(T)\setminus\mathcal E$.
Since $T$ is $r$-uniform and $r\ge 3$, no edge in
$E(T)\setminus\mathcal E$ is contained in $N_T[e]$.
Hence $E(T[N_T[e]])=\mathcal E$,
and $|f\cap N_T[e]|\le 1$ for any
$f\in E(T)\setminus E(T[N_T[e]])$.
\end{proof}

\begin{proposition}\label{prop:defect}
Let $P$ be a maximum edge-neighbourhood packing of $T\in\Tr$,
$p=|P|=\iota(T)$ and $\ell=|V(T)\setminus\bigcup_{e\in P}N_T[e]|$. Then
\begin{equation}\label{eq:defect}
n(T)-(2r-1)p=\ell+(r-1)\sum_{e\in P}(a(e)-1).
\end{equation}
Thus, $\iota(T)\le\lfloor n(T)/(2r-1)\rfloor$.
\end{proposition}

\begin{proof}
Since $N_T[e]\cap N_T[f]=\emptyset$ for all distinct
$e,f\in P$, by Lemma~\ref{lem:neighbourhood}, we have
\[
\begin{aligned}
n(T)
&=\ell+\sum_{e\in P}|N_T[e]|\\
&=\ell+\sum_{e\in P}\bigl(r+(r-1)a(e)\bigr)\\
&=(2r-1)p+\ell+(r-1)\sum_{e\in P}(a(e)-1).
\end{aligned}
\]

Since $\ell\ge 0$ and $a(e)\ge 1$ for every $e\in P$, $(2r-1)p\le n(T)$.
Hence, $\iota(T)=p\le\lfloor n(T)/(2r-1)\rfloor$.
\end{proof}

\section{The extremal hypertrees}
We first define the marked hypertrees needed for the
extremal characterization.
Let $\alpha_1,\ldots,\alpha_r$ be nonnegative integers with
$k=\sum_{j=1}^r\alpha_j\ge 1$, and $C=\{c_1,\ldots,c_r\}$ be an edge.
For each $j\in\{1,\ldots,r\}$, add $\alpha_j$ edges,
each consisting of $c_j$ and $r-1$ new vertices.
All new vertices are distinct.
Let $M_r(\alpha_1,\ldots,\alpha_r)$ be the resulting
marked hypertree with distinguished edge $C$.
We call $C$ the \emph{central edge} and the vertices in
$V(M_r(\alpha_1,\ldots,\alpha_r))\setminus C$
the \emph{attachment vertices}.
Then $n\bigl(M_r(\alpha_1,\ldots,\alpha_r)\bigr)=r+k(r-1)$.

\begin{figure}
\centering
\includegraphics{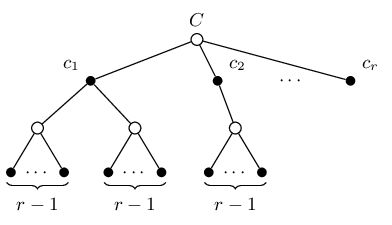}
\caption{The incidence tree of $M_r(2,1)$.}
\label{fig:marked}
\end{figure}

Permuting $\alpha_1,\ldots,\alpha_r$ gives an isomorphic
marked hypertree. For $1\le t<r$ and nonnegative integers
$\alpha_1,\ldots,\alpha_t$ with
$\sum_{j=1}^{t}\alpha_j\ge 1$, let
\[
M_r(\alpha_1,\ldots,\alpha_t)
:=
M_r(\alpha_1,\ldots,\alpha_t,
\underbrace{0,\ldots,0}_{r-t}).
\]
Define
\[
\begin{aligned}
\M_0(r)&=\{M_r(1)\},\\
\M_1(r)&=\{M_r(2),M_r(1,1)\},\\
\M_2(r)&=\{M_r(3),M_r(2,1),M_r(1,1,1)\}.
\end{aligned}
\]
These families contain six marked isomorphism types.
We have $m(M)=d+2$ and $n(M)=2r-1+(r-1)d$
for every $d\in\{0,1,2\}$ and $M\in\M_d(r)$.
For each marked hypertree $M$ with central edge $C$,
only the vertices in $V(M)\setminus C$ may be used
for further attachments in the construction below.
Figure~\ref{fig:marked} shows the incidence graph of $M_r(2,1)$.

\begin{construction}\label{con:extremal}
Let $p$ and $s$ be integers with $p\ge 1$ and
$0\le s\le 2r-2$ such that $(2r-1)p+s\equiv 1\pmod{r-1}$.
Choose $d_1,\ldots,d_p\in\{0,1,2\}$ and an integer
$\ell\ge 0$ satisfying
\begin{equation}\label{eq:construction}
\ell+(r-1)\sum_{i=1}^{p}d_i=s.
\end{equation}

Let $M_1,\ldots,M_p$ be pairwise vertex-disjoint marked
hypertrees with $M_i\in\M_{d_i}(r)$ for $1\le i\le p$,
and $C_i$ be the central edge of $M_i$.
Let $U$ be a vertex set satisfying $|U|=\ell$ and
$U\cap\bigcup_{i=1}^{p}V(M_i)=\emptyset$.
Let $z_1,\ldots,z_p$ be distinct new vertices, and
$F$ be an $r$-uniform hypertree with $V(F)=\{z_1,\ldots,z_p\}\cup U$.
When $p+\ell=1$, $F$ consists of a single vertex and no edges.

For each $f\in E(F)$ and each $i$ with $z_i\in f$,
choose a vertex $x_{i,f}\in V(M_i)\setminus C_i$.
We allow $x_{i,f}=x_{i,g}$ for distinct
$f,g\in E(F)$ with $z_i\in f\cap g$.
For each $f\in E(F)$, replace each $z_i\in f$ with
the chosen attachment vertex $x_{i,f}$ of $M_i$
and retain the vertices in $f\cap U$.
Let $f'$ be the resulting edge, that is,
\[
f'
=
(f\cap U)\cup
\{x_{i,f}:1\le i\le p,\ z_i\in f\}.
\]
Let $T$ be the hypergraph with
$V(T)=U\cup\bigcup_{i=1}^{p}V(M_i)$ and
$E(T)=\bigcup_{i=1}^{p}E(M_i)\cup\{f':f\in E(F)\}$.
Let $\Ep$ be the family of unmarked hypergraphs obtained
by this construction, and call $F$ a \emph{skeleton} of $T$.
\end{construction}

The skeleton always exists.
Indeed, the congruence condition and~\eqref{eq:construction}
imply $p+\ell\equiv p+s\equiv 1\pmod{r-1}$.
Starting with a single vertex, successively add
$(p+\ell-1)/(r-1)$ edges, each consisting of one existing
vertex and $r-1$ new vertices.
The resulting hypergraph is an $r$-uniform hypertree
of order $p+\ell$. When $p+\ell=1$, no edge is added.

Each $f'$ has size $r$ and contains at most one
vertex of each $M_i$.
Replacing these vertices by the corresponding $z_i$
recovers $f$, so the edges defining $T$ are distinct.
The incidence graph $B(T)$ is obtained from the tree
$B(F)$ by replacing each $z_i$ with the tree $B(M_i)$
and reattaching the incident edges at the chosen vertices.
Thus $B(T)$ is a tree.
\begin{theorem}\label{thm:extremal}
Let $T\in\Tr$ have order $n$, $p=\lfloor n/(2r-1)\rfloor$ and
$s=n-(2r-1)p$. Then $\iota(T)=\left\lfloor\frac{n}{2r-1}\right\rfloor$
if and only if $T\in\Ep$.
\end{theorem}

\begin{proof}
Suppose that $T\in\Ep$. Let $F,M_1,\ldots,M_p$ and $U$ be as in
Construction~\ref{con:extremal}.
Let $C_i$ be the central edge of $M_i$ for $1\le i\le p$.
By construction, $f'\cap C_i=\emptyset$
for any $f\in E(F)$ and $1\le i\le p$.
It follows that $N_T[C_i]=V(M_i)~(1\le i\le p)$.
Since the sets $V(M_1),\ldots,V(M_p)$ are pairwise
disjoint, $\{C_1,\ldots,C_p\}$ is an edge-neighbourhood
packing of $T$.
Therefore, $\iota(T)\ge p$ by Theorem~\ref{thm:packing}.
By~\eqref{eq:construction}, we have
\[
n=\sum_{i=1}^{p}\bigl(2r-1+(r-1)d_i\bigr)+\ell
=(2r-1)p+s.
\]
Since $0\le s\le 2r-2$, by Proposition~\ref{prop:defect},
\[
p\le\iota(T)\le
\left\lfloor\frac{n}{2r-1}\right\rfloor=p.
\]

Conversely, suppose that $\iota(T)=p$.
There exists a maximum edge-neighbourhood packing
$P=\{e_1,\ldots,e_p\}$ of $T$ by
Theorem~\ref{thm:packing}.
For $1\le i\le p$, let $M_i=T[N_T[e_i]]$ with
central edge $e_i$ and $d_i=a(e_i)-1$.
Let $U=V(T)\setminus\bigcup_{i=1}^{p}N_T[e_i]$,
$\ell=|U|$.
Each $d_i$ is a nonnegative integer by
Lemma~\ref{lem:neighbourhood}.
By Proposition~\ref{prop:defect}, we have
\begin{equation}\label{eq:reverse}
\ell+(r-1)\sum_{i=1}^{p}d_i=s\le 2r-2.
\end{equation}
Since $\ell\ge 0$, $\sum_{i=1}^{p}d_i\le 2$.
Therefore, $d_i\in\{0,1,2\}$ for any $i$.

For each $i$, $M_i$ consists of $e_i$ and $d_i+1$
other edges by Lemma~\ref{lem:neighbourhood}.
Since $B(T)$ is a tree, each of these edges contains
one vertex of $e_i$ and $r-1$ vertices outside $e_i$,
and the sets of vertices outside $e_i$ are pairwise
disjoint.
Since $d_i+1\le 3$, we have $M_i\in\M_{d_i}(r)$
for $1\le i\le p$ by the definitions of the marked families.
Let $\mathcal R=E(T)\setminus\bigcup_{i=1}^{p}E(M_i)$.
By Lemma~\ref{lem:neighbourhood}, we have $|f\cap V(M_i)|\le 1$ and $f\cap e_i=\emptyset$ for any $f\in\mathcal R$ and $1\le i\le p$.

For distinct $i,j$, we have
$V(M_i)\cap V(M_j)=\emptyset$ and
$E(M_i)\cap E(M_j)=\emptyset$.
Thus the graphs $B(M_i)$ are pairwise
vertex-disjoint subtrees of $B(T)$.
Let $B'$ be the tree obtained from $B(T)$ by contracting
each subtree $B(M_i)$ to a vertex $z_i$.
For each $f\in\mathcal R$, let
\[
f'
=
(f\cap U)\cup
\{z_i:1\le i\le p,\ f\cap V(M_i)\ne\emptyset\}.
\]
The neighbours in $B'$ of the vertex corresponding
to $f$ are exactly the elements of $f'$.
Since $|f\cap V(M_i)|\le 1$ for any $i$, we have
\[
|f'|
=
|f\cap U|+\sum_{i=1}^{p}|f\cap V(M_i)|
=
|f|
=
r.
\]
If $f'=g'$ for distinct $f,g\in\mathcal R$,
then the two vertices corresponding to $f$ and $g$
and any two distinct vertices of $f'$ form a
$4$-cycle in $B'$.
This contradicts the fact that $B'$ is a tree.

Let $F$ be the hypergraph with $V(F)=\{z_1,\ldots,z_p\}\cup U$
and $E(F)=\{f':f\in\mathcal R\}$.
Then $B(F)\cong B'$, so $F$ is an $r$-uniform hypertree.
If $E(F)=\emptyset$, then $F$ consists of a single
vertex, as allowed in Construction~\ref{con:extremal}.

For each $f\in\mathcal R$ and $i$ with $z_i\in f'$,
let $x_{i,f'}$ be the unique vertex of $f\cap V(M_i)$.
Then $x_{i,f'}\in V(M_i)\setminus e_i$ and
$f=(f'\cap U)\cup\{x_{i,f'}:1\le i\le p,\ z_i\in f'\}$.
Thus, $T$ can be obtained by Construction~\ref{con:extremal}
with these choices and the parameters satisfying~\eqref{eq:reverse}.
Therefore, $T\in\Ep$.
\end{proof}

\begin{proof}[Proof of Theorem~\ref{thm:bound}]
The inequality follows from Proposition~\ref{prop:defect}.
For an admissible order $n$, write $n=(2r-1)p+s$, where $0\le s\le2r-2$.
In Construction~\ref{con:extremal}, choose $d_i=0$ for every $i$ and $\ell=s$.
Since $p+s\equiv n\equiv1\pmod{r-1}$, a skeleton exists.
Theorem~\ref{thm:extremal} gives a hypertree of order $n$ with the isolation number $p$.
\end{proof}
Table~\ref{tab:types} lists the possible numbers of marked
hypertrees from each of the families $\M_0(r)$, $\M_1(r)$
and $\M_2(r)$, together with the corresponding values of $\ell$.
Here $p\ge 1$ and $p+s\equiv 1\pmod{r-1}$.
The last row requires $p\ge r$, since $p\ge 2$ and
$p\equiv 1\pmod{r-1}$ in this case.

\begin{table}[htbp]
\centering
\caption{Possible distributions of the marked hypertrees
in Construction~\ref{con:extremal}.}
\label{tab:types}
\begin{tabular}{ccccc}
\hline
$s$ & $\M_0(r)$ & $\M_1(r)$ & $\M_2(r)$ & $\ell$\\
\hline
$0\le s\le 2r-2$
  & $p$ & $0$ & $0$ & $s$\\
$r-1\le s\le 2r-2$
  & $p-1$ & $1$ & $0$ & $s-r+1$\\
$s=2r-2$
  & $p-1$ & $0$ & $1$ & $0$\\
$s=2r-2$
  & $p-2$ & $2$ & $0$ & $0$\\
\hline
\end{tabular}
\end{table}
\begin{corollary}\label{cor:unrounded}
Let $T\in\Tr$.
Then $\iota(T)=n(T)/(2r-1)$ if and only if
$T\in\mathcal{E}^{(r)}_{p,0}$ for some positive integer
$p$ with $p\equiv 1\pmod{r-1}$.
Equivalently, $T$ can be obtained by
Construction~\ref{con:extremal} with $U=\emptyset$
and $M_i\in\M_0(r)$ for $1\le i\le p$.

For an integer $m\ge 2$, there exists a member of $\Tr$
with $m$ edges attaining this equality if and only if
$m=2+(2r-1)j$ for some integer $j\ge 0$.
\end{corollary}

\begin{proof}
Let $p=\lfloor n(T)/(2r-1)\rfloor$ and
$s=n(T)-(2r-1)p$.
Since $\iota(T)$ is an integer, the equality
$\iota(T)=\frac{n(T)}{(2r-1)}$ holds if and only if $s=0$
and $T\in\Ep$ by Theorem~\ref{thm:extremal}.
In Construction~\ref{con:extremal}, $s=0$ if and only if
$\ell=0$ and $d_i=0$ for every $i$
by~\eqref{eq:construction}.
These conditions are equivalent to $U=\emptyset$ and
$M_i\in\M_0(r)$ for every $i$.

For such a construction, let $F$ be its skeleton
and $j=m(F)$.
Since $\ell=0$, the skeleton $F$ has $p$ vertices, so
$p=1+(r-1)j$.
Thus $p\equiv 1\pmod{r-1}$.
Each $M_i$ has two edges, so $m(T)=2p+j=2+(2r-1)j$.

Conversely, let $j\ge 0$ be an integer and
$p=1+(r-1)j$.
Choose an $r$-uniform hypertree $F$ of order $p$.
Let $H$ be obtained by Construction~\ref{con:extremal}
with the skeleton $F$, $U=\emptyset$ and
$M_i\in\M_0(r)$ for $1\le i\le p$.
Then $m(H)=2+(2r-1)j$ and $n(H)=(2r-1)p$.
By Theorem~\ref{thm:extremal}, $\iota(H)=p=\frac{n(H)}{(2r-1)}$.
\end{proof}
Let $r=4$.
Choose pairwise vertex-disjoint marked hypertrees
$M_1,M_2\in\M_1(4)$ and $M_3,M_4\in\M_0(4)$.
Let $U=\emptyset$ and $F$ be the hypertree with the
single edge $\{z_1,z_2,z_3,z_4\}$.
Let $T$ be obtained by Construction~\ref{con:extremal}
with these choices.
Here $p=4$, $\ell=0$ and $s=6$.
We have $n(T)=2\cdot10+2\cdot7=34$ and
$m(T)=2\cdot3+2\cdot2+1=11$.
Hence, by Theorem~\ref{thm:extremal},
\[
\iota(T)=4=\left\lfloor\frac{34}{7}\right\rfloor
<\frac{34}{7}.
\]

\begin{remark}\label{rem:nonunique}
Every extremal hypertree can be obtained by
Construction~\ref{con:extremal}.
The same hypertree may be obtained from different
choices in the construction.
\end{remark}

\section{Triangle block graphs}
A set $S\subseteq V(G)$ is a \emph{cycle isolating set}
of $G$ if $G-N_G[S]$ is a forest.
The \emph{cycle isolation number} $\iota_c(G)$ is the
minimum cardinality of a cycle isolating set of $G$.
Borg~\cite{BorgCycles} proved the sharp bound
$\iota_c(G)\le \frac{n}{4}$ for every connected graph $G$ of
order $n$ with $G\not\cong K_3$.

A vertex $v$ of a graph $G$ is a \emph{cut vertex}
if $G-v$ has more components than $G$.
A \emph{block} of a graph is a maximal connected
subgraph with no cut vertex.
We obtain the following result for connected graphs
whose blocks are all triangles.

\begin{corollary}\label{cor:triangles}
Let $G$ be a connected graph of order $n$ with at least
two blocks, all of which are triangles.
Let $T$ be the hypergraph with $V(T)=V(G)$ and
$E(T)=\{V(B):B\text{ is a block of }G\}$.
Then $\iota_c(G)\le\left\lfloor\frac n5\right\rfloor$.
The bound is attained for every odd integer $n\ge5$.
Equality holds if and only if
$T\in\mathcal{E}^{(3)}_{p,s}$, where
$p=\lfloor n/5\rfloor$ and $s=n-5p$.
\end{corollary}
\begin{proof}
Since $G$ is connected, $B(T)$ is connected.
Distinct blocks of $G$ have at most one vertex in common,
so $B(T)$ contains no $4$-cycle.
Suppose that $B(T)$ contains a cycle of length at least $6$.
Then $G$ contains a cycle with edges from different blocks.
This contradicts the fact that every cycle of $G$
is contained in a single block.
Thus, $B(T)$ is a tree.
Since every block of $G$ is a triangle and $G$ has
at least two blocks, $T\in\mathcal{T}_3$.

Let $S\subseteq V(G)$.
Two distinct vertices of $G$ are adjacent if and only
if they belong to a common edge of $T$.
Therefore, $N_G[S]=N_T[S]$.
Every cycle of $G$ is contained in a block, and every
block of $G$ is a triangle.
Thus $G-N_G[S]$ is a forest if and only if
$e\cap N_T[S]\ne\emptyset$
for every $e\in E(T)$.
It follows that $S$ is a cycle isolating set of $G$
if and only if $S$ is an isolating set of $T$.
By Theorem~\ref{thm:bound}, we have
\[
\iota_c(G)=\iota(T)\le
\left\lfloor\frac n5\right\rfloor.
\]
Equality holds if and only if
$T\in\mathcal{E}^{(3)}_{p,s}$
by Theorem~\ref{thm:extremal}.

For any odd integer $n\ge5$, there exists a
hypertree $H\in\mathcal{T}_3$ of order $n$ satisfying
$\iota(H)=\lfloor n/5\rfloor$
by Theorem~\ref{thm:bound}.
Let $G_H$ be the graph obtained from $H$ by replacing
each edge with a triangle on the same vertex set.
Since $B(H)$ is a tree, $G_H$ is connected and its
blocks are exactly these triangles.
The graph $G_H$ has $n$ vertices and at least two blocks.
By the equality established above,
\[
\iota_c(G_H)=\iota(H)
=\left\lfloor\frac n5\right\rfloor.
\]
\end{proof}
The condition that every block is a triangle excludes bridges.
The bound in Corollary~\ref{cor:triangles} need not hold
if bridges are allowed.
For example, let $G$ be obtained from two vertex-disjoint
triangles $C_1$ and $C_2$ by adding a path of length $3$
with one end in each triangle and two new internal vertices.
We have $n(G)=8$ and $N_G[V(C_1)]\cap N_G[V(C_2)]=\emptyset$.
Thus every cycle isolating set of $G$ contains at least
two vertices.
The two ends of the added path form a cycle isolating set, so
\[
\iota_c(G)=2>\left\lfloor\frac85\right\rfloor.
\]

Let $T\in\Tr$, and let $G$ be the graph obtained by
replacing each edge of $T$ with a copy of $K_r$
on the same vertex set, where $K_r$ denotes the
complete graph of order $r$.
The blocks of $G$ are exactly these copies of $K_r$.
A set $S\subseteq V(G)$ is a \emph{$K_r$-isolating set}
of $G$ if $G-N_G[S]$ contains no copy of $K_r$.
For any $S\subseteq V(T)$, we have $N_G[S]=N_T[S]$.
Thus $S$ is an isolating set of $T$ if and only if
$S$ is a $K_r$-isolating set of $G$.
When $r=3$, this is equivalent to $G-N_G[S]$ being
a forest, as in Corollary~\ref{cor:triangles}.
For $r\ge4$, the graph $G-N_G[S]$ may contain a triangle
even when it contains no copy of $K_r$.

\end{document}